\documentclass[11p,leqno]{amsart}
\usepackage{amsmath}
\usepackage{amsfonts}
\usepackage{amssymb}
\usepackage{graphicx}
\usepackage{color}
\newtheorem{theorem}{Theorem}[section]
\newtheorem*{theorem*}{Theorem}
\newtheorem{lemma}{Lemma}[section]
\newtheorem{corollary}[theorem]{Corollary}
\newtheorem{proposition}{Proposition}[section]

\newtheorem{remark}[theorem]{Remark}

\def\qed{\hfill $\square$}

\def\qed{\hfill $\square$}

\newcommand{\RR}{\mathbb{R}}

\def\heat{\lf(\frac{\p}{\p t}-\Delta\ri)}

\def\Ric{\text{Ric}}
\def\lf{\left}
\def\ri{\right}

\def\p{\partial}

\def\R{\RR}

\def\Sph{\mathbb S}

\def\tpsi{\tilde{\psi}}

\def\id{\operatorname{id}}
\def\Ric{\operatorname{Ric}}
\def\Rm{\operatorname{Rm}}

\newcommand{\PP}{\mathcal{P}}               
\newcommand{\Ptwo}{\PP_2(M)}            

\numberwithin{equation}{section}

\begin{document}

\title[Sharp logarithmic Sobolev]{\bf Sharp logarithmic Sobolev inequalities on gradient solitons and applications}
\author{Jos\'e A. Carrillo and Lei Ni}

\begin{abstract}
We show that gradient shrinking, expanding  or steady  Ricci
solitons have potentials leading to suitable reference probability
measures on the manifold.  For shrinking solitons, as well as
expanding soltions with nonnegative Ricci curvature, these
reference measures satisfy sharp logarithmic Sobolev inequalities
with lower bounds characterized by the geometry of the manifold.
The geometric invariant appearing in the sharp lower bound is
shown to be  nonnegative. We also  characterize the expanders when
such invariant  is zero. In the proof various useful volume growth
estimates are also established for gradient shrinking and
expanding solitons. In particular, we prove that the {\it
asymptotic volume ratio} of any gradient shrinking soliton with
nonnegative Ricci curvature must be zero.
\end{abstract}

\maketitle

\section{Introduction}

A complete Riemannian manifold  $(M, g)$ is called a gradient
shrinking  soliton (shrinker) if there exists a (smooth) function
$f$, such that its Hessian $f_{ij}$ satisfies
\begin{equation}\label{soliton-mother}
R_{ij}+f_{ij} -\frac{1}{2}g_{ij}=0.
\end{equation}
Here $R_{ij}$ denotes the Ricci curvature. As shown in Theorem 4.1
of \cite{CLN}, associated to the metric and the {\it potential
function} $f$, there exists a family of metrics $g(\eta)$, a
solution to Ricci flow
$$
\frac{\partial}{\partial \eta} g(\eta) =-2\Ric(g(\eta)),
$$
with the property that $g(0)=g$, the original metric,   and a
family of diffeomorphisms $\phi(\eta)$, which is generated by the
vector field  $X=\frac{1}{\tau} \nabla f$, such that $\phi(0)=\id$
and $g(\eta)=\tau(\eta) \phi^*(\eta) g$ with $\tau(\eta)=1-\eta$,
as well as $f(x, \eta)=\phi^*(\eta) f(x)$. Namely there exists a
self-similar (shrinking) family of metrics which is a solution to
the Ricci flow. The metric $g(\eta)$ and $f(\eta)$, sometimes also
written as $g^\tau$ and $f^\tau$, or simply $g$ and $f$ when the
meaning is clear, satisfy that
\begin{equation}\label{soliton-f}
R_{ij}+f_{ij} -\frac{1}{2\tau}g_{ij}=0.
\end{equation}
We shall denote by $S(x)$ the scalar curvature and by
$d\Gamma_\tau$ the volume element of $g^\tau$.

Gradient shrinking solitons arise as the singularity models of
Ricci flow. The more interesting cases are the noncompact ones.
Trivial examples include the Euclidean space $\R^n$ and the
cylinders $\Sph^k\times \R^{n-k}$ for $k\ge 2$. Non-trivial
noncompact examples can be found in, for example \cite{FIK}. There
is also a more recent construction of solitons with symmetry in
\cite{dw}. The main result of this paper is  the following
theorem, which generalizes the sharp Logarithmic Sobolev
Inequality (LSI) of the Euclidean space $\R^n$ \cite{Gr}. This was
referred as Stam-Gross logarithmic Sobolev inequality in
\cite{Villani,V08}, where one can also find detailed historic
accounts and more complete references.

\begin{theorem}\label{main1}
Assume that $(M, g, f)$ is a gradient shrinking soliton, then:
\begin{enumerate}
\item[i)] The potential $e^{-f}$ is integrable on $M$ and it can
be normalized as
\begin{equation}\label{norm}
\frac{1}{(4\pi \tau)^{n/2}}\int_M e^{-f}d\Gamma_\tau =1.
\end{equation}

\item[ii)] LSI inequality: There exists a geometric invariant
$\mu_s$, under isometries, which depends only on the value of $f$
and $S$ at the minimum point of $f$ and is independent of the
parameter $\tau$, such that
$$
\int_M \left\{\tau\left[|\nabla \psi|^2_{\tau}+S(\cdot,
\tau)\right]+\psi -n\right\}\rho\, d\Gamma_\tau \ge -\mu_s,
$$
for any $\tau>0$ and any nonnegative compactly supported smooth
function $\rho=\frac{e^{-\psi}}{(4\pi \tau)^{n/2}}$ with unit
integral on $M$. Moreover, the geometric invariant $\mu_s$ in the
above inequality is sharp.

\item[iii)] If $|\Rm|\le A$ for some $A>0$, then $\mu_s\ge0$.
\end{enumerate}
\end{theorem}

We refer to Section 2 for notations and the definition of the
invariant $\mu_s$. Let us observe that the expression in the LSI
makes sense wherever $\rho=0$ (hence $\psi=-\infty$) since $\rho
\, \psi =0$ there. This can be seen more easily if the integrand
is expressed in terms of $\rho$,
$$
\int_M \left\{4 \tau |\nabla \sqrt{ \rho}|^2_{\tau}+\rho S(\cdot,
\tau) - \rho \ln \rho -n\rho -\frac{n}2 \ln (4\pi
\tau)\rho\right\}\, d\Gamma_\tau \ge -\mu_s.
$$
Note also that for the Gaussian soliton, namely $(M, g, f)=(\R^n,
g_{\operatorname{can}}, \frac{1}{4}|x|^2)$, we get $\mu_s=0$. The
proof of the inequality uses the Bakry-Emery criterion \cite{BE}
for LSI's, as obtained from the so-called HWI inequalities derived
by Villani and coauthors in different settings \cite{V08}, see
Section 3 and references therein. Here, the main difficulty
resides in showing the necessary bounds on the potential to
normalize it as in \eqref{norm} and thus, being able to apply
these inequalities, which is done in Section 2. We should point
out that Perelman has claimed \cite[Remark 3.2]{P} that when a LSI
holds on a shrinking solitons, the sharp form can be justified
using his entropy formula. This mainly applies to compact
shrinkers since for the noncompact case, even a weak form of LSI
is not known. For the compact shrinkers, our approach supplies a
different argument. An immediate consequence of the theorem is the
strong non-collapsing of the gradient shrinking solitons. In the
case that $M$ has bounded nonnegative Ricci curvature the sharp
LSI of Theorem \ref{main1} implies LSI's for all scales, not
necessarily with sharp constants though, hence the non-collapsing
at all scales. Let us finally comment that the sharpness in the
third statement of our main theorem means that taking the density
$\rho$ to be the normalized potential in the first statement, then
the minimum $\mu_s$ is attained, see Section 4.

In our analysis of the gradient shrinking solitons we also prove
the following result.

\begin{corollary}\label{coro-main} Any non-flat gradient shrinking soliton with nonnegative
Ricci curvature must have zero asymptotic volume ratio.
\end{corollary}

This  is also done in Section 2. This result, in the case of
gradient shrinking solitons, generalizes a previous  result of
Perelman \cite{P} on ancient solutions with bounded nonnegative
curvature operator. The result of Perelman \cite[Proposition
11.4]{P} draws the same conclusion for any ancient solutions with
bounded nonnegative curvature operator. Let us remark that
Perelman also assumes the non-collapsing condition, which however
is not needed in the proof. Gradient shrinking solitons are
special ancient solutions. However our result is more general than
that of Perelman for the following reasons: it cannot be derived
from Perelman's since we do not assume that the curvature operator
is nonnegative nor bounded only Ricci curvature must be
nonnegative. On the other hand, the result of Perelman can be
derived out of the result above on gradient steady solitons  via
the {\it asymptotic solitons} \cite[Proposition 11.2]{P}.

Gradient steady/expanding solitons (expanders) arise also in the
singularity analysis  of Ricci flow \cite{H93}. The technique
employed here yields similar sharp geometric inequalities for
 expanding solitons as well.

\begin{theorem} Assume that $(M, g, f)$ is a gradient expanding soliton with
$\Ric \ge 0$. Then for any nonnegative
$\rho(x)=\frac{e^{-\psi(x)}}{(4\pi)^{n/2}}$ with $\int_M \rho(x)\,
d\Gamma(x)=1$, we have that
\begin{equation}\label{lsi-gesi}
\int_M \left(|\nabla \psi|^2-3S+\psi-n\right)\rho\, d\Gamma \ge
-\mu_e.
\end{equation}
Here $\mu_e$ is a geometric constant  depending only on the value
of $f$ and $S$ at the minimum point of $f$. The
inequality is sharp for such $\mu_e$. Moreover $\mu_e\ge 0$ with
equality if and only if $(M, g)$ is isometric to $\R^n$.
\end{theorem}

We refer the readers to Section 5 for  notations involved in the above theorem.
An equivalent expression of the integrand appeared in a
recent interesting preprint of Cao and Hamilton \cite{Cao-Hamilton} on
pointwise differential estimates of Li-Yau-Hamilton type.

For  expanding solitons,  we also obtain a  volume estimate, which
generalizes a recent result of Hamilton \cite{H05}, see also
\cite{CLN}, asserting that {\it the asymptotic volume ratio of
gradient expanding solitons with bounded positive Ricci curvature
must be positive}. The following is one of our statements.

\begin{corollary} Assume that $(M, g, f)$ is a gradient expanding
soliton with $S(x)\ge -\beta$ for some constant $\beta\ge 0$. Then
for any $o\in M$ and $r\ge r_0$
$$
V(o, r)\ge V(o, r_0)\left(\frac{r+a}{r_0+a}\right)^{n-2\beta}
$$
with $a=2\sqrt{f(o)+\mu_e +\beta}$.
\end{corollary}

The above mentioned Hamilton's result follows from the above
statement applying to the case $\beta=0$. For general $\beta$, the
growth rate in our estimate  is sharp as shown by examples. This
is proved in Section 5. A similar/independent result can also be
found in a recent preprint \cite{ChT}.

For gradient steady soliton, since one can not expect that the LSI
holds in general in viewing of Hamilton's `cigar', we obtain a
sharp weighted Poincar\'e inequality instead. The proof is
relatively easy, without appealing the above mentioned theory
involving the Bakry-Emery criterion, and is done in Section 6.

The part  $\mu_s,\mu_e\ge0$ of the main theorems is finally proved
in Section 7. This is motivated by the Zamolodchikov's c-theorem
of the re-normalization group flow \cite{Zam}.  In view of the
entropy monotonicity formula of Perelman, its connection with the
LSI, as well as the fact that gradient shrinking solitons arise as
the singularity models (at least for the  cases that  the blow-up
has  nonnegative curvature operator), this result can be viewed as
an analogue of Zamolodchikov's c-theorem
 for the re-normalization group flow. The proof makes
use a Li-Yau-Hamilton type inequality of Perelman \cite{P} and the
 entropy formula/monotonicity for the linear heat equation of \cite{N1}.

By the work of Dolbeault and Del Pino \cite{DoPi} (also Toscani
and the first author \cite{CT})  the sharp form of
Sobolev-Gagliardo-Nirenberg inequalities on $\R^n$ are related to
the nonlinear Fokker-Planck dynamics of porous medium/fast
diffusion type equations. It is interesting to find out if such
relation can lead to sharp inequalities on solitons along a
similar line of argument of this paper.

\section{Geometric estimates on gradient shrinking solitons}

We shall follow the notations of the introduction being our first
objective to show the integrability on the potential function for
solitons allowing for the normalization \eqref{norm}. The
following equations are simple consequences of the soliton
equation (\ref{soliton-f}):
\begin{eqnarray}
S+\Delta f -\frac{n}{2\tau}&=&0,\label{soliton1}\\
S+|\nabla f|^2 -\frac{f}{\tau} &=&\frac{\mu_s(\tau)}{\tau}. \label{soliton2}
\end{eqnarray}
where $\mu_s(\tau)$ is a constant that will be eventually chosen
by the normalization of the potential as in \eqref{norm}. Here $S$
is the scalar curvature. See, for example \cite{H93} or
\cite{CLN}, for a proof. The equations (\ref{soliton1}) and
(\ref{soliton2}) imply that
\begin{equation}\label{soliton3}
2\Delta f -|\nabla f|^2 +S+\frac{f-n}{\tau} =-\frac{\mu_s(\tau)}{\tau}.
\end{equation}

The lemma below implies that the integral involved in the
normalization \eqref{norm}, as well as other integrals involved
later in the proof of Theorem \ref{main1}, are finite.

\begin{lemma}\label{lembehaviorinfty}
 Let $r(x)$ be the distance function to a fixed point
$o\in M$ with respect to $g(\eta)$ metric. Then there exist
constant $C_1$ and $C_2$ such that
\begin{equation}\label{help4}
f(x)\ge \frac{1}{4\tau} (r(x)-C_1)^2
\end{equation}
for $r(x)\ge C_2$ and
\begin{equation}\label{help5}
f(x)\le \frac{1}{4\tau}\left( r(x)+C_1\right)^2, \quad \quad  |\nabla f|(x) \le \frac{1}{2\tau}(r(x)+C_1)
\end{equation}
for $r(x)\ge C_2$.
\end{lemma}
\noindent {\sl Proof.} First we observe that $S\ge 0$ by a
gradient estimate argument of Chen \cite{Ch} (see also the
appendix of \cite{Ta}). The estimate (\ref{help4}) then follows
verbatim from \cite[pages 655--656]{FMZ}.

Now (\ref{soliton2}) and $S\geq 0$ imply that
\begin{equation}\label{aux}
|\nabla f|\le \sqrt{\frac{f+\mu_s(\tau)}{\tau}}.
\end{equation}
The first estimate in (\ref{help5}) follows easily from this by
integrating $f+\mu_s(\tau)$ along minimizing geodesics from $o$,
see also the proof of Proposition 5.1. The second estimate in
(\ref{help5}) follows from the first one via \eqref{aux}. \qed

\begin{corollary}\label{corshrsol}
Let $(M, g, f)$ be a gradient shrinking soltion. Then the minimum of $f$ can
be achieved somewhere, say $o\in M$. Hence $f(o)$ and $S(o)$ are
fixed for different $\eta$ and the constant $\mu_s(\tau)$ in
\eqref{soliton2} is a constant independent of $\tau$. Therefore
$\mu_s\doteqdot \mu_s(\tau)$ is an invariant of the soliton.
Moreover,
\begin{equation}\label{int-help}
\int_M \left(|\Delta f|+|\nabla f|^2 +|f|+S\right)e^{-f}\, d\Gamma_\tau <\infty.
\end{equation}
\end{corollary}

\noindent {\sl Proof.} The first part of the corollary is evident
since $o$ is the fixed point of $\phi(\eta)$. The change of $S$
(from the shrinking) is compensated by the factor $\tau$. For the
second part, observe that by \cite[page 7]{WW}, we have that
$$
e^{-f}J(r, \theta) \le A_1 e^{a_2 r-\frac{1}{2\tau}r^2}
$$
for some positive constants $A_1, a_2$ independent of $r$. Here
$J(r, \theta)$ is the area element of the sphere $\partial
B_o(r)$. Namely, $\operatorname{Area}(\partial
B_o(r))=\int_{\Sph^{n-1}}J(r, \theta)\,d \theta.$ Notice that
(\ref{soliton1}) and (\ref{soliton2}), along with (\ref{help4})
and (\ref{help5})  effectively imply that
\begin{eqnarray}
0\le S &\le& \frac{1}{4\tau^2} (r(x)+C_1)^2, \label{scalarupper}\\
-\frac{n}{2}\le-\Delta f &\le& \frac{1}{4\tau^2} (r(x)+C_1)^2. \label{laplace}
\end{eqnarray}
Since $\int_M (\cdot) \, d\Gamma_\tau=\int_0^\infty
\int_{\Sph^{n-1}}(\cdot)J(r, \theta)\, d\theta\, dr$, the
finiteness of the integral in (\ref{int-help}) follows easily.
\qed

Note that by Theorem 4.1 of \cite{CLN}, $\frac{e^{-f}}{(4\pi
\tau)^{n/2}}$ satisfies the conjugate heat equation
$$
\left(\frac{\partial}{\partial \tau} -\Delta +S\right)
\left(\frac{e^{-f}}{(4\pi \tau)^{n/2}}\right)=0.
$$
Hence the total mass of $e^{-f}$, namely the normalization
(\ref{norm}) is preserved along the evolution. In other words, if
$$
\int_M \frac{e^{-f}}{(4\pi\tau)^{\frac{n}{2}}}\,d\Gamma_\tau=1
$$
holds at $\tau=1$ (which corresponds to $\eta=0$), it holds for
all $\tau>0$. Also note that $\mu_s(M, g)=\mu_s(M', g')$ if $(M,
g)$ is isometric to $(M', g')$  by the virtue of \cite[Lemma
1.2]{Naber}. Equivalently, the invariant $\mu_s(M, g)$ is
independent of the choice of the potential function $f$ since the
difference of two potential functions is either a constant or a
linear function, since they have the same Hessian. In the first
case, the normalization make the constant zero. For the second
case, namely the difference of the two potential functions is a
linear function, then the manifold $M$  splits off a line. Some
simple calculation also show that the normalization would make the
constants $\mu_s$ in (\ref{soliton2}) identical for the two
different potential functions. In fact, if the difference of two
potential functions $h \doteqdot f_1-f_2$ (assuming $\tau=1$
without the loss of the generality), is a linear function of $\R$
and $M= \R\times M_1$, using the soliton equation (\ref{soliton1})
one can write $f_k(x, y)=\frac{1}{4}x^2 +b_k x+c_k(y)$ for
$k=1,2$. Here we denote the coordinate of $\R$ by $x$ and the
coordinate of $M'$ by $y$. Since $h(x)=ax+b$ for constant $a$ and
$b$. Hence we have $c_1(y)-c_2(y)=c$. Now if $\int_M
e^{-f_1}=\int_M e^{-f_2}$, by simple direct calculation we have
that $b_1^2=c+b_2^2$. Direct calculation shows that
$$
\left(|\nabla f_1|^2+S-f_1\right)-\left(|\nabla f_2|^2+S-f_2\right)=b_1^2-b_2^2 -(c_1(y)-c_2(y))=0.
$$

The next result concerns the behavior of the volume $V(o, r)$ of
balls $B(o, r)$, especially as $r\to \infty$. We start with the
easier case of the  Ricci curvature being bounded.

\begin{corollary}\label{cor:volnoopt}
Let $(M, g)$ be a nonflat gradient shrinking soliton with
$\Ric\ge0$. Assume further that the scalar curvature $S(x)\le C_1$
for some $C_1>0$. Then, there exists a constant $\delta=\delta(M,
f)>0$ with the property that for any $o\in M$, there exists
$a=a(M, f, C_1)$ such that for any $r\ge r_0\ge a$
\begin{equation}\label{vol-est1}
V(o, r)\le V(o, r_0)\left(\frac{r-a}{r_0-a}\right)^{n-\delta}.
\end{equation}
\end{corollary}
\noindent {\sl Proof.} Without the loss of generality we may
assume that $\tau=1$. By \cite[Proposition 1.1]{N}, we have
$\delta=\delta(M, f)>0$ such that $S\ge \frac{\delta}{2}$. On the
other hand, by \cite[Section 8]{P}, see also the proof of
\cite[Proposition 1.1]{N}, for any minimizing geodesic joining $o$
to $x=\gamma(s_0)$ with $s_0\ge 2$ and $r_0> 0$ such that
$s_0-r_0\ge 1$, we have that
\begin{equation}\label{est3}
\int_0^{s_0-r_0} \Ric(\gamma', \gamma')\, ds \le C_4(M)+\frac{n-1}{r_0}.
\end{equation}
This implies, again by the argument in the proof of
\cite[Proposition 1]{N}, that
$$
\frac{\partial f}{\partial r}(x) \ge \frac{r(x)}{2}-C_6(M, f, o).
$$
Now integration by parts on equation (\ref{soliton1}) over $B(o,
r)$ yields that
\begin{eqnarray*}
\frac{n-\delta}{2}V(o, r)&\ge & \int_{B(o, r)} \left(\frac{n}{2}-S\right) \, d\Gamma \\
&=&\int_{\partial B(o, r)}\frac{\partial f}{\partial r}\, dA \\
&\ge& A(o, r)\left(\frac{r}{2}-C_6\right).
\end{eqnarray*}
Here $A(o, r)$ is the surface area of $\partial B(o, r)$. The
result follows from integrating the above estimate on $[r_0, r]$. \qed

\begin{remark}
Being Corollary {\rm \ref{cor:volnoopt}} proved under no
restriction on the boundedness of the Ricci curvature, it might be
used to prove \cite[Proposition 11.4]{P}. This result concludes
that {\it any nonflat ancient solution of Ricci flow  with bounded
nonnegative curvature operator must have the asymptotic volume
ratio $\lim_{r\to \infty} \frac{V(o, r)}{r^n}=0$.} In fact one can
derive \cite[Proposition 11.4]{P} by contradiction: Assume the
claim is false, one obtains an {\it asymptotic soliton} by
\cite[Proposition 11.2]{P}, which is nonflat and that has the
maximum volume growth. On one hand, we now may use the volume
comparison theorem, to get that $V(o, r)/r^n$ is always bounded
above the corresponding ratio of the Euclidean space. On the other
hand, it is easy to show that the asymptotic volume ratio for an
ancient solution with bounded nonnegative curvature is monotone
non-increasing in $t$. Hence the asymptotic soliton must has
positive asymptotic volume ratio. This is a contradiction with the
estimate \eqref{vol-est1}.
\end{remark}

With some extra effort, we can indeed prove such desired volume
estimate without assuming the Ricci curvature upper bound.

\begin{proposition}\label{gss-vol-est}
Let $(M, g, f)$ be a nonflat gradient shrinking soliton with
$\Ric\ge0$. Then
$$
\lim_{r\to \infty} \frac{V(o, r)}{r^n}=0.
$$
Here $V(o, r)$ is the volume of $B(o, r)$.
\end{proposition}
\noindent {\sl Proof.} We can reduce ourselves to the case
$\tau=1$ without loss of generality. For simplicity, after
translation we may assume that the potential function $f$
satisfies $ |\nabla f|^2 +S\le f. $ It is also more convenient to
work with sub-level sets of $f$. Let us consider the sets
$$
F_r \doteqdot \left\{x\in M\,|\, 2\sqrt{f(x)}\le r\right\}
$$
and $\widetilde{V}(r)=\operatorname{Vol}(F_r)$. Assume that the
conclusion is not true, then $\liminf_{r\to \infty}
r^{-n}V(o,r)\ge \eta
>0$ for some $\eta>0$. This clearly implies that $\liminf_{r\to
\infty} r^{-n}\widetilde{V}(r)\ge\eta'>0$.

On the other hand, following \cite{CZh,Mu} we consider the
function
$$
\chi(r)=\int_{F_r} S\, d\Gamma.
$$
Using $S\ge \delta>0$ for some $\delta>0$, which is ensured by
\cite[Proposition 1.1]{N}, we have that for any small $\epsilon>0$
\begin{eqnarray*}
(n-\delta')\widetilde{V}(r)+2(1-\epsilon)\chi(r) &=&
\int_{F_r} 2\left(\frac{n}{2}-\epsilon \delta+(1-\epsilon)S\right)\, d\Gamma \\
&\ge&\int_{F_r} 2(\frac{n}{2}-S)\, d\Gamma \\
&=& \int_{F_r} 2\Delta f\, d\Gamma\\
 &=&r\widetilde{V}'(r) -\frac{4\chi'(r)}{r}.
\end{eqnarray*}
Here $\delta'=2\epsilon \delta$ and in the last line we used the
computation in (4) of \cite{Mu}. Integrating the above estimate as
in \cite{Mu,CZh}, we arrive at
$$
\frac{\widetilde{V}(r)}{r^{n-\delta'}} -
\frac{\widetilde{V}(r_0)}{r_0^{n-\delta'}}\le
\frac{4\chi(r)}{r^{n-\delta'+2}}
$$
for $r\ge r_0\ge 8\sqrt{n+2}$. Now using that $2\chi(r)\le
n\widetilde{V}(r)$, we have that the right hand side above tends
to zero as $r\to \infty$. This induces  that $\limsup_{r\to
\infty}r^{-n} \widetilde{V}(r)=0$, which is a contradiction. \qed

We should remark that there exists a proof to Perelman's result by
Hamilton via his singularity analysis of ancient solutions. The
interested reader can find the details of Hamilton's argument in
\cite{CLN}.
 It is interesting to find out if  Proposition \ref{gss-vol-est} can be shown for any ancient solutions with nonnegative Ricci curvature.

\section{Optimal Transport and LSIs}

In this section, we will work with Riemannian manifolds $(M,g)$
endowed with a reference probability measure $e^{-V} d\Gamma$
where the potential $V\in C^2(M)$ verifies a curvature-dimension
bound of the type $C(K,\infty)$ with $K\in \RR$, i.e.,
$$
R_{ij}+V_{ij} \geq K g_{ij} .
$$
Here $d\Gamma$ is the volume measure associated to $(M,g)$. This
section is devoted to collect several results present in the
literature \cite{V08}. A Riemannian
manifold in this section refers to a smooth, complete connected
finite-dimensional Riemannian manifold distinct from a point,
equipped with a smooth metric tensor. Let us assume that the
reference measure is normalized by
$$
\int_M e^{-V} d\Gamma = 1.
$$
Consider the positive solution $\rho$ to the Fokker-Planck
equation
\begin{equation}\label{eq1.1}
\frac{\partial \rho}{\partial t} -\operatorname{div} \left(\rho
\nabla(\log \rho + V) \right)=0.
\end{equation}
Let $\xi=\log \rho + V$. It is easy to see that
$$
\heat \xi =\langle \nabla \xi, \nabla \log \rho\rangle.
$$
Let us define the Boltzmann relative entropy functional, called
also Nash entropy, as
$$
H_{V} (\rho)\doteqdot\int_M \rho \xi\, d\Gamma.
$$
We have immediately the following dissipation of the Boltzmann
relative entropy functional,
\begin{equation}
\frac{d}{dt}H_{V}(\rho(t))= -\int_M |\nabla \xi|^2 \rho\, d\Gamma
\doteqdot - I_{V}(\rho(t)), \label{1std}
\end{equation}
where computations are made for smooth, fast-decaying at infinity
for non-compact manifolds, solutions on $M$. This computation show
us that these two quantities, the relative Boltzmann entropy
$H_{V} (\rho)$ and the {\it relative Fisher information}
$I_V(\rho)$ are intimately related at least for solutions of
\eqref{eq1.1}. However, as it was discovered in the case of
$\R^n$,  and in the case of a manifold in \cite{BE,Arnold} for
linear diffusions or in \cite{CT,Otto,DoPi} for nonlinear
diffusions, this relation is really through functional
inequalities, see also \cite {Otto-Villani}.

Related to these functionals, there is another quantity that is
involved in these inequalities: the Euclidean Wasserstein distance
between any two probability measures $\nu_0$, $\nu_1$ on the
manifold $M$, i.e.,
\begin{equation}\label{W2}
W_2(\nu_0,\nu_1) \doteqdot \inf \left \{ \int_{M\times M}
r^2(x,y)\,d\theta(x,y); \ \theta \in \Theta(\nu_0,\nu_1)\right
\}^{1/2};
\end{equation}
where $\Theta(\nu_0,\nu_1)$ is the set of probability measures on
$M\times M$ having marginals $\nu_0$ and $\nu_1$, $r(x, y)$ is the
Riemannian distance between $x$ and $y$. This distance is well
defined for probability measures $\nu_0$ and $\nu_1$ with second
moment bounded, $\Ptwo$, and metrizes the weak convergence of
measures in the sense of \cite[Definition 6.7, Theorem 6.8]{V08}.
The expression ``second moment bounded" refers to the fact that
the squared distance function $r^2(x)$ is integrable against the
measures $\nu_0$ and $\nu_1$. It worths to mention that the
curvature-dimension bound $C(K,\infty)$ with $K>0$ implies that
the second moment of the reference measure (actually, all moments)
$e^{-V}$ is bounded, see \cite[Theorem 18.11]{V08}.

Recently, several authors \cite{S,LV} based on early works
\cite{Mc,CMS}, see \cite[Chapter 17]{V08} for a whole account of
the history, have characterized curvature-dimension bounds in
terms of the displacement convexity of the Boltzmann relative
entropy functional. The notion of displacement convexity refers to
convexity along pathes of minimal transport distance $W_2$ in the
set of probability measures $\Ptwo$. An expression of the
convexity of these functionals are the so called HWI inequalities,
named in this way since they involved the three functionals $H_{V}
(\rho)$, $I_V(\rho)$ and $W_2$. In the following, we will work
with measures absolutely continuous against volume measure and we
identify the measures with their densities for notational
convenience. The main results we need are the following:

\begin{theorem} \label{thmlsi} {\rm (\cite[Corollary 20.13]{V08} and \cite{BE})}
Let M be a Riemannian manifold equipped with a reference measure
$e^{-V} d\Gamma$ where the potential $V\in C^2(M)$ verifies a
curvature-dimension bound of the type $C(K,\infty)$ with $K\in\R$.
Then, for any given $\nu\in\Ptwo$ absolutely continuous with
respect to volume measure $d\Gamma$ with density $\rho$, it holds
the HWI inequality:
$$
H_V(\rho)\leq W_2(\rho,e^{-V})\sqrt{I_V(\rho)} -
\frac{K}{2} W_2(\rho,e^{-V})^2 .
$$
As a consequence, we have that whenever $K>0$, the following LSI
follows
$$
H_V(\rho)\leq \frac{1}{2K} I_V(\rho).
$$
\end{theorem}

The HWI inequalities were originally introduced in
\cite{Otto-Villani} and used in other models in nonlinear PDEs in
\cite{CMV}. Later, they were generalized to compact manifolds in
\cite{LV} and in this generality in \cite{V08}. To see that the
LSI inequality follows from the HWI inequality it suffices to
consider the right-hand side of the HWI inequality as a function
of $W_2$ and maximize that function.

Let us remark that some proofs of the LSI inequality use the
Fokker-Planck dynamics \eqref{eq1.1}, called the Bakry-Emery
stragegy, but the referred functional proof through the HWI
inequalities allows to overcome discussions on integrability
issues and the decay at infinity for non-compact manifolds of
solutions to \eqref{eq1.1}. In fact, a direct application of the
LSI on \eqref{1std} gives the exponential decay of the Boltzmann
relative entropy functional for solutions of \eqref{eq1.1} with
initial density in $\Ptwo$ in case $C(K,\infty)$ with $K>0$ holds,
i.e., given a solution $\rho(t)$ of \eqref{eq1.1} then
$$
H_V(\rho(t)) \leq H_V(\rho(0))\, e^{-2Kt} \qquad \mbox{for all }
t\geq 0.
$$
Nevertheless, let us remind the reader that assuming all
integrability and behavior at the infinity are met for all integration
by parts below, we can obtain the evolution of the relative Fisher
information (see also \cite{BE,Arnold,Villani,V08} for these
computations). To take the time derivative of $I_V(\rho(t))$ note
the Bochner type formula
\begin{equation}\label{bochner1}
\heat |\nabla \xi|^2 =-2\xi_{ij}^2+2\langle \nabla (\langle \nabla
\xi, \nabla \log \rho\rangle), \nabla \xi\rangle -2R_{ij}\xi_i
\xi_j.
\end{equation}
Using the above formula we have that
\begin{eqnarray*}
\frac{d}{dt} I_V(\rho(t)) &=& \int_M (\Delta |\nabla \xi|^2)\rho
+|\nabla \xi|^2 \operatorname{div}(\nabla \rho
+\rho\nabla V)\, d\Gamma\\
&\quad & +\int_M \left(-2\xi_{ij}^2+2\langle \nabla (\langle
\nabla \xi, \nabla \log \rho\rangle), \nabla \xi\rangle
-2R_{ij}\xi_i \xi_j\right)\rho\, d\Gamma.
\end{eqnarray*}
Since
\begin{eqnarray*}
\int_M \langle \nabla (\langle \nabla \xi, \nabla \log
\rho\rangle), \nabla \xi\rangle\rho\, d\Gamma &=& \int_M
\langle \nabla (|\nabla \xi|^2-\langle \nabla V, \nabla \xi\rangle), \nabla \xi\rangle\rho\, d\Gamma\\
&=&\int_M
\langle \langle \nabla |\nabla \xi|^2, \nabla \rho\rangle +\langle \nabla |\nabla \xi|^2, \nabla V\rangle \rho\, d\Gamma\\
&\quad& -\int_M\langle \nabla \langle \nabla V, \nabla \xi\rangle,
\nabla \xi\rangle\rho\, d\Gamma
\end{eqnarray*}
we arrive at
\begin{eqnarray}
\frac{d}{dt} I_V(\rho(t))&=&\int_M \left(-2\xi_{ij}^2
-2R_{ij}\xi_i \xi_j\right)\rho\, d\Gamma\nonumber\\
&\quad& +\int_M \langle \nabla |\nabla \xi|^2, \nabla V\rangle
\rho
-2\langle \nabla \langle \nabla V, \nabla \xi\rangle, \nabla \xi\rangle\rho\, d\Gamma\nonumber \\
&=&\int_M \left(-2\xi_{ij}^2 -2(R_{ij}+V_{ij})\xi_i
\xi_j\right)\rho\, d\Gamma \label{2ndd}.
\end{eqnarray}
As a consequence, due to the curvature dimension bound
$C(K,\infty)$, we have
$$
\frac{d}{dt} I_V(\rho(t)) \leq -2K \int_M |\nabla \xi|^2\rho\, d\Gamma,
$$
and thus,
$$
I_V(\rho(t)) \leq I_V(\rho(0))\, e^{-2Kt} \qquad \mbox{for all }
t\geq 0.
$$

\section{Main Result and Applications}

Now, let us come back to the precise situation we have, the case
of a shrinking soliton,  and prove the main Theorem \ref{main1}.
Let us define the potential $V=f+\frac{n}2 \log(4\pi \tau)$ for
the fixed time slice of the shrinking Riemannian manifold soliton
$(M, g)$ at time $\tau$. Lemma \ref{lembehaviorinfty} and
Corollary \ref{corshrsol} implies that $e^{-V}$ is a well defined
probability measure. Moreover, we deduce from the soliton
definition \eqref{soliton-f} that this reference measure verifies
the $C(\frac{1}{2\tau},\infty)$ condition. Therefore, Theorem
\ref{thmlsi} implies that for any probability density of the form
$$
\rho(x)=\frac{e^{-\psi(x)}}{(4\pi\tau)^{\frac{n}{2}}}
$$
with second moment bounded, we get the LSI
$$
H_V(\rho)\leq \tau I_V(\rho).
$$
Using now the soliton equation \eqref{soliton3}, we deduce:
\begin{eqnarray*}
I_V(\rho)&=&\int_M \left(|\nabla \psi|^2 \rho +2\langle \nabla f,
\nabla \rho\rangle +|\nabla f|^2 \rho\right)\, d\Gamma\\
&=& \int_M \left[|\nabla \psi|^2 \rho +(-2\Delta f+|\nabla
f|^2)\rho\right]\,
d\Gamma\\
&=& \int_M \left[|\nabla \psi|^2+S+\frac{f+\mu_s
-n}{\tau}\right]\rho \, d\Gamma.
\end{eqnarray*}
Thus, the LSI inequality is equivalent to
\begin{equation}\label{lsi}
\int_M \left[\tau(|\nabla \psi|^2+S)+\psi -n\right]\rho\, d\Gamma
\ge -\mu_s,
\end{equation}
for all densities $\rho$ with bounded second moment for the
shrinking soliton, with $\mu_s$ characterized by Corollary
\ref{corshrsol}.

Now recall Perelman's entropy functional
$$
\mathcal{W}(g^\tau, u, \tau)\doteqdot \int_M \left[\tau(|\nabla
\psi|^2+S)+\psi-n\right] u\, d\Gamma_\tau
$$
is defined for $u=\frac{e^{-\psi}}{(4\pi \tau)^{n/2}}$ with
$\int_M u\, d\Gamma_\tau=1$. Theorem \ref{main1} implies that for
$(M, g^\tau)$, $W(g^\tau, u, \tau)\ge -\mu_s$. Namely Perelman's
$\mu$-invariant
$$
\mu(g^\tau, \tau)\doteqdot \inf_{\int_M u=1}\mathcal{W}(g^\tau,
u,\tau)
$$
is bounded from below by $-\mu_s$. From (\ref{soliton3}) it is
easy to see that
$$
\tau(2\Delta f-|\nabla f|^2+S)+f-n=-\mu_s.
$$
Hence $u=\frac{e^{-f}}{(4\pi \tau)^{n/2}}$ is the minimizer for
Perelman's $\mu(g, \tau)$, cf. \cite[Remark 3.2]{P}. This shows
that the inequality of  Theorem \ref{main1} is sharp. Summarizing,
we have that

\begin{corollary}
Let $(M, g, f)$ be a gradient shrinking soliton satisfying \eqref{soliton-mother}. Then
  $$
  \mu(g, 1)=-\mu_s.
  $$
\end{corollary}

\begin{remark}
When $f=\operatorname{constant}$,  $(M, g)$ is a Einstein manifold
with $\Ric_M =\frac{1}{2}g_M$. In this case we obtain a
log-Sobolev inequality for $S=\frac{n}{2}$ and
$$
\mu_s=\frac{n}{2}-\log (V(M))+\frac{n}{2}\log (4\pi)
$$
where $V(M)$ is the volume of $(M, g_M)$. The $\mu$-invariant  was
computed in \cite{CHI} for many examples of four manifolds.

When $M=\R^n$ with $f=\frac{1}{4} |x|^2$, direct calculation shows
that $\mu_s=0$. Hence the  classical logarithmic Sobolev
inequality of Stam-Gross is a special case.
\end{remark}

Recall here that a solution of Ricci flow is called $\kappa$
non-collapsed, if for any $(x_0, t_0)$ and $r\ge 0$, such that on
$P(x_0, t_0, r)=B_{g(t_0)}(x_0, r)\times [t_0-r^2, t_0]$,
$|Rm|(x,t)\le r^{-2}$, then $V_{g(t_0)}(x_0, r)\ge \kappa r^n$.
Here $V_{g(t_0)}(x_0, r)$ is the volume of $B_{g(t_0)}(x_0, r)$
with respect to $g(t_0)$. Perelman \cite[Theorem \ref{main1}]{P}
implies the following volume non-collapsing result for gradient
shrinking solitons.

\begin{corollary}
Let $(M, g, f)$ be a gradient shrinking soliton satisfying
\eqref{soliton-mother}. Then there exists a
$\kappa=\kappa(\mu_s)>0$ such that if in a ball $B(x_0, 1)$,
$|\Ric|\le 1$, then $V(x_0, 1)\ge \kappa$. In particular, if the
Ricci curvature is bounded on $M$ which is noncompact, then $M$
has at least linear volume growth.
\end{corollary}

\noindent {\sl Proof.} Follows from Theorem \ref{main1} and
Section 4 of \cite{P}. See also \cite{Chowetc, topping}. \qed

In \cite{Naber} there is a related result asserting the
$\kappa$-noncollapsing of gradient shrinking solitions with
bounded curvature, in the sense  defined right above the
corollary. The conclusion in above corollary appears stronger
since it only requires global lower bound on the scalar curvature
and the local bound of the Ricci curvature over the ball, for a
fixed time-slice only.

When $\Ric(M, g_{\tau=1}) \ge 0$ and is bounded, one can
derive the logarithmic Sobolev inequality for all scales. This is
done in the following two propositions.

\begin{proposition}[Scale $> 1$] \label{as-lsi1}
Let $(M , g)$ be a gradient shrinking soliton satisfying
\eqref{soliton-mother}. Assume that $\Ric\ge 0$. Then, there
exists positive $\delta=\delta(M)<1$ such that for any $\sigma>
1$,
$$
\int_M \left[\sigma (|\nabla \tpsi|^2+S)+\tpsi
-n\right]\frac{e^{-\tpsi}}{(4\pi \sigma)^{\frac{n}{2}}}\, d\Gamma
\ge -\mu_s+\frac{n}{2}-\delta -\frac{n}{2}\log\left(
\frac{n}{2\delta}\right)
$$
for any $\tpsi$ satisfying that $\int_M e^{-\tilde\psi}/{(4\pi
\sigma)^{\frac{n}{2}}}\,d\Gamma=1$.
\end{proposition}

\noindent {\sl Proof.} Clearly only the nonflat case worths the
proof (since the flat one is isometric to $\R^n$). By
\cite[Proposition 1.1]{N}, for a nonflat gradient shrinking
soliton, there exists $\delta=\delta(M, f)>0$ such that $S(x)\ge
\delta$ for any $x\in M$. Let $\psi=\tpsi +\frac{n}{2}\log
\sigma$. Then it is easy to see that
\begin{eqnarray*}
\int_M \left[\sigma (|\nabla \tpsi|^2+S)+\tpsi -n\right]\frac{e^{-\tpsi}}{(4\pi \sigma)^{\frac{n}{2}}}\, d\Gamma
&=& \int_M \left(|\nabla \psi|^2+S+\psi -n\right)\frac{e^{-\psi}}{(4\pi)^{\frac{n}{2}}}\, d\Gamma\\
&\, &+(\sigma -1) \int_M (|\nabla \psi|^2+S)\frac{e^{-\psi}}{(4\pi)^{\frac{n}{2}}}\, d\Gamma -\frac{n}{2}\log \sigma\\
&\ge & -\mu_s +\delta(\sigma-1) -\frac{n}{2}\log \sigma,
\end{eqnarray*}
where we have used Theorem \ref{main1} in the last estimate. Since
$\delta(\sigma -1)-\frac{n}{2}\log \sigma\ge \frac{n}{2}-\delta
-\frac{n}{2}\log( \frac{n}{2\delta})$, the claimed result
follows.\qed

From the proof, the following corollary is evident, observing that
$S\ge 0$ for shrinking solitons, which is clear from \cite{Ch},
see also the appendix of  \cite{Ta}.

\begin{corollary}
Let $(M , g)$ be a gradient shrinking soliton satisfying
\eqref{soliton-mother}.  Then for any $\sigma>1$,
$$
\int_M \left[\sigma (|\nabla \tpsi|^2+S)+\tpsi
-n\right]\frac{e^{-\tpsi}}{(4\pi \sigma)^{\frac{n}{2}}}\, d\Gamma
\ge -\mu_s-\frac{n}{2}\log \sigma.
$$
\end{corollary}

\begin{proposition}[Scale $< 1$] \label{as-lsi2}
Assume that $0\le \Ric\le A$. Then for any $0\le \sigma\le 1$,
$$
\int_M \left[\sigma (|\nabla \tpsi|^2+S)+\tpsi
-n\right]\frac{e^{-\tpsi}}{(4\pi \sigma)^{\frac{n}{2}}}\, d\Gamma
\ge -\mu_s-nA
$$
for any $\tpsi$ satisfying that $\int_M e^{-\tilde\psi}/{(4\pi
\sigma)^{\frac{n}{2}}}\,d\Gamma=1$.
\end{proposition}

\noindent {\sl Proof.} Define
$$
\mu_0(g, \sigma)\doteqdot \inf_{\int_M u_0=1} \int_M \left(\sigma
|\nabla \tpsi|^2+\tpsi-n\right)u_0 \,d\Gamma
$$
with $u_0= e^{-\tpsi}/{(4\pi \sigma)^{\frac{n}{2}}}$. Theorem
\ref{main1} implies that $\mu_0(g, 1) \ge -\mu_s -nA.$ Now for any
$u_0$ which is compactly supported, let $u(x,t)$ be the heat
equation solution with $u(x, 0)=u_0$. Then  by the entropy
monotonicity result in \cite{N1}, for $\sigma\le 1$,
\begin{eqnarray*}
\int_M \left(\sigma |\nabla \tpsi|^2+\tpsi-n\right)\frac{e^{-\tpsi}}{(4\pi \sigma)^{\frac{n}{2}}}
&\ge& \int_M \left(|\nabla \varphi|^2+\varphi-n\right)u(y, 1-\sigma)\, d\Gamma(y)\\
&\ge&\mu_0(g, 1)
\end{eqnarray*}
where $u(y, 1-\sigma)=e^{-\varphi(y)}/{(4\pi
(1-\sigma))^{\frac{n}{2}}}$. This implies the claimed result. \qed

The above two propositions imply that $\nu(g)>-\infty$, see
section 7 for a definition, hence the strong
$\kappa$-non-collapsing result for gradient shrinking solitons
with bounded and nonnegative Ricci curvature as in \cite{P} (see
also \cite{topping} and \cite{Chowetc}). For the general case
without assuming $\Ric\ge 0$, one can still obtain a logarithmic
Sobolev for scales less than one, see Section 7.

\section{Expanding solitons}

Recall that $(M, g)$ is called a gradient expanding soliton if
there exists $f$ such that
\begin{equation}\label{soliton-mother2}
R_{ij}+\frac{1}{2}g_{ij}=f_{ij}.
\end{equation}
It is easy to show that
\begin{eqnarray}
\Delta f&=&S+\frac{n}{2}\label{ges-1}\\
S+|\nabla f|^2-f&=&\mu_e\label{ges-2}
\end{eqnarray}
for some constant $\mu_e$. As before we will eventually choose
$\mu_e$ by the normalizing condition $\int_M
e^{-f}/{(4\pi)^{n/2}}\, d\Gamma =1$. This will make $\mu_e$ a
geometric invariant of $(M, g)$. T

Our first concern is about the behavior of the volume of balls
$B(o, r)$ in  $M$ for any given $o\in M$. Along this direction, Hamilton \cite{H05}
proved the following result:

\begin{theorem} Let $(M, g)$ be a gradient expanding soliton  has bounded
nonnegative Ricci curvature.  Then $(M, g)$ has maximum volume
growth. Namely
$$
\liminf_{r\to \infty}\frac{V(o, r)}{r^n} >0.
$$
\end{theorem}

For the exposition of this result please see \cite[Proposition
9.46]{CLN}. Let us remark that the assumption of uniform
boundedness of the Ricci curvature is used in the proof to bound
$\int_{\gamma} \Ric(\gamma' \gamma')$ as in Section 2. Here, the
limit always exists due to the Bishop-Gromov volume comparison.
The limit of the quotient is called the  {\it asymptotic volume
ratio}. This compares sharply with the gradient shrinking solitons
(cf. Proposition \ref{gss-vol-est}) and a result of Perelman
\cite{P} asserting that {\it any non-flat ancient solution with
bounded nonnegative curvature operator has zero asymptotic volume
ratio}. The result below is a generalization of the above result
of Hamilton.

\begin{proposition}\label{hamilton-avr}
Let $(M, g, f)$ be an gradient expanding soliton.

\noindent {\rm (1)}  If $S(x)\ge 0$ for any $x\in M$, without
assuming any curvature bound, then for any $o\in M$, $r\ge r_0$.
$$
V(o, r)\ge V(o, r_0)\left(\frac{r+a}{r_0+a}\right)^n
$$
 with $a=2\sqrt{f(o)+\mu_e}$.

\noindent  {\rm (2)} Assume that $S(x)\ge -\beta$ for some
constant $\beta>0$. Then for any $o\in M$ and $r\ge r_0$
$$
V(o, r)\ge V(o, r_0)\left(\frac{r+a}{r_0+a}\right)^{n-2\beta}
$$
with $a=2\sqrt{f(o)+\mu_e +\beta}$.
\end{proposition}

\noindent {\sl Proof.}  In the case (1),  from the assumption and
(\ref{ges-2}) we have that $f+\mu_e\ge 0$. Consider any minimizing
geodesic $\gamma(s)$ from $o\in M$ a fixed point of $M$. Then
(\ref{ges-2}) implies  that for any $s$
$$
\left|\frac{d}{ds} f(\gamma(s))\right|^2\le f+\mu_e.
$$
This  implies, by the ODE comparison, that
$$
\left(2\sqrt{f+\mu_e}\right)(\gamma(s))\le s+a
$$
where $a=2\sqrt{f(o)+\mu_e}$, which then implies that
\begin{equation}\label{ges-est2}
\left|\frac{\partial f}{\partial r}\right|(\gamma(s)) \le
\frac{s}{2}+\frac{a}{2}.
\end{equation}

Now we integrate (\ref{ges-1}) on $B(o, r)$ and have that
\begin{eqnarray*}
\frac{n}{2}V(o, r)&\le& \frac{n}{2} V(o, r)+\int_{B(o, r)}S\, d\Gamma \\
 &=& \int_{B(o, r)} \Delta f\, d\Gamma\\
&\le&\int_{\partial B(o, r)} \left|\frac{\partial f}{\partial
r}\right|(y)\, dA(y).
\end{eqnarray*}
Using  (\ref{ges-est2}) we have that
$$
\frac{n}{2}V(o, r)\le A(o, r)(\frac{r}{2}+\frac{a}{2}).
$$
The result follows by dividing the both side of the above by $V(o,
r)$ and then integrating the resulting estimate on the interval
$[r_0, r]$. The proof for the case (2) is similar. \qed

\begin{remark}
The estimates in both cases have the sharp power. To see this
consider $M =N^k \times \R^{n-k}$ where $N$ is a compact Einstein
manifold with $\Ric_N =-\frac{1}{2}g_N$, $\R^{n-k}$ is the
Gaussian expanding soliton.
\end{remark}

Now we derive the LSI for the expanders. To make sure that the
integral $\int_M e^{-f}\, d\Gamma$ is finite we have to make an
assumption that there exists some $\epsilon >0$,
\begin{equation}\label{ges-basic1}
f_{ij}=\frac{1}{2}g_{ij}+R_{ij}\ge \epsilon g_{ij}.
\end{equation}
Under this assumption, it is easy to see that
$$
f(x)\ge \frac{\epsilon}{4}r^2(x) -C
$$
for some $C=C(M, f)$. Since $R_{ij}\ge -\frac{1}{2}g_{ij}$,  the
volume $V(o, r)\le \exp(A(r+1))$ for some $A=A(n)$. This together
with the lower estimate above ensures that the integral $\int_M
e^{-f}\, d\Gamma$ is finite, see also \cite{WW}. Notice that under
our assumption (\ref{ges-basic1}), as in the proof of Proposition
\ref{hamilton-avr} we have that
$$
f(x)\le \left(\frac{r(x)}{2}+b \right)^2
$$
for some $b=b(M, f)$. This ensures the finiteness of the
integral
$$\int_M \left(|\nabla f|^2+|\Delta f|+|S|\right)\frac{e^{-f}}{(4\pi)^{n/2}}\, d\Gamma.$$
Note that (\ref{ges-1}) and (\ref{ges-2}) implies that
\begin{equation}\label{ges-3}
2\Delta f -|\nabla f|^2 -3S +f-n =-\mu_e.
\end{equation}
Integrating (\ref{ges-3}), we have that
$$
\int_M \left(|\nabla f|^2-3S
+f-n\right)\frac{e^{-f}}{(4\pi)^{n/2}}\, d\Gamma =-\mu_e.
$$
It is clear that assumption \eqref{ges-basic1} is trivially
satisfied for the case that $M$ has non-negative Ricci curvature.

Assume in the rest of this section that $\Ric \ge 0$, let us
define the potential $V=f-\frac{n}2 \log(4\pi)$. Previous
arguments imply that the reference measure $e^{-V}$ is a well
defined probability measure. Moreover, we deduce from the soliton
definition \eqref{soliton-mother2} and being $\Ric \ge 0$ that
this reference measure verifies the $C(\frac12,\infty)$ condition.
Therefore, Theorem \ref{thmlsi}, together with a similar calculation as before,  implies the following LSI
inequality.

\begin{theorem}\label{ges-main1}
Assume that $(M, g, f)$ is gradient expanding soliton with $\Ric
\ge 0$. Then for any $\rho(x)=e^{-\psi(x)}/{(4\pi)^{n/2}}$ with
$\int_M \rho(x)\, d\Gamma(x)=1$, we have that
\begin{equation}\label{lsi-ges}
\int_M \left(|\nabla \psi|^2-3S+\psi-n\right)\rho\, d\Gamma \ge
-\mu_e.
\end{equation}
\end{theorem}

Here $\mu_e$, as before, is a geometric invariant (in the sense of
Section 2), which is the same for two isometric metrics. One can
write in the dynamic form by considering the family of metrics
$g(\tau)$ (in this case with $g(1)$ being the original metric, and
$0<\tau<\infty$) generated by the diffeomorphisms, as for the
shrinking solitons case described in the introduction. Since it is
the same inequality by re-scaling we omit its full statement. Note
that in the left hand side of (\ref{lsi-ges}) an equivalent
integrand is
$$
\tau\left(2\Delta \psi-|\nabla \psi|^2 -3S\right)+\psi-n.
$$
This expression also showed itself up in a differential Harnack or
Li-Yau-Hamilton type calculation, in a recent preprint of Cao and
Hamilton \cite{Cao-Hamilton}, where however  the nonnegativity of
the curvature operator is required. It is certainly interesting to
explore the connections between the LSI here and the
Li-Yau-Hamilton type estimate for Ricci flow solution.

\begin{corollary}
Let $(M, g)$ be an expanding soliton as in Theorem
\ref{ges-main1}. Then $M$ is diffeomorphic  to $\R^n$.
\end{corollary}
\noindent {\sl Proof.} First it is easy to see that $M$ is of
finite topological type. This follows from the observation that
$f$ is a proper function and has no critical point outside a
compact subset \cite{FMZ},  since for any $x\in M$ and $\gamma(s)$
a minimizing geodesic jointing $o\in M $, a fixed point, to $x$,
with $f(\gamma(0))=o$ and $f(\gamma(s_0))=x$
\begin{eqnarray*}
f'(\gamma(s_0))&=&f'(\gamma(0))+\int_0^{s_0} f''(\gamma(s))\, ds\\
&\ge& \frac{s_0}{2}+f'(\gamma(0)).
\end{eqnarray*}
The conclusion follows from the uniqueness of the critical point
along with the strict convexity of $f$. \qed

\section{Gradient steady solitons}

Now we consider the  gradient steady solitons. Recall that a
 gradient steady soliton $(M, g)$ has a potential function $f$
satisfying that
\begin{equation}\label{sgs1}
R_{ij}=f_{ij}.
\end{equation}
It was shown in \cite{H93} that
\begin{equation}\label{sgs2}
|\nabla f|^2 +S=\lambda
\end{equation}
for some $\lambda$.
Similar as before there is a solution to Ricci flow $g(\tau)$
associated with  the gradient steady soliton $(M, g, f)$ \cite{CLN}. We first
need the following lemma to ensure the finiteness of $\int_M
e^{-f}\, d\Gamma$ and other integrals later involved, under some
geometric assumptions.

\begin{lemma}\label{gs-help1}
Let $(M, g, f)$ be a gradient steady soliton. Assume that there
exists a point $o\in M$ such that $S(o)=\max_M S$ and either
$\Ric(x)>0$ for all $x\in M$,  or $\Ric\ge 0$ and
$$
\limsup_{x\to \infty} S(x) <\max_M S.
$$
Then $o$ is a minimum of $f$ and there exists $\delta>0$ and $C=C(M, f)$ so that
\begin{equation}\label{sgs-lb1}
f(x)\ge \delta r(x) -C.
\end{equation}
Here $r(x)$ is the distance function to $o$. In particular, $M$ is
diffeomorphic to  $\R^n$ in the case $\Ric>0$ and of finite
topological type in the case $\Ric\ge 0$.
\end{lemma}

\noindent {\sl Proof.} For the first case, it was shown in
\cite[Theorem 20.1]{H93} that $o$ is the unique minimum of $f$.
Note that the argument there actually requires $\Ric>0$ even
though it was not stated; it is also necessary, as shown by easy
examples. Note that for any geodesic $\gamma(s)$ from $o$, we have
that
$$
\frac{d^2}{ds^2}(f(\gamma(s))=\Ric(\gamma', \gamma')>0.
$$
Hence we have for any $s_0>0$, $\frac{d}{ds}(f(\gamma(s_0)))>0$.
Then $f(\gamma(s)) \ge
\frac{d}{ds}(f(\gamma(s_0)))(s-s_0)+f(\gamma(s_0))$, which implies
the desired lower estimate.

For the second case, the assumption already excludes the Ricci
flat situation, on which clearly (\ref{sgs-lb1}) fails for $f$
being a  constant. We first claim that under the assumption on the
behavior of $S$ at the infinity, $S(o)=\lambda$. Suppose it is not
true, then $\max_M S<\lambda$ and $|\nabla f|^2\ge \lambda-\max_M
S$. Let $\sigma(u)$ be an integral curve of $\nabla f$ passing $o$
with $\sigma(0)=o$. Direct calculation shows that
$\frac{d}{du}\left(|\nabla f|^2(\sigma(u))\right)=2\Ric(\nabla f,
\nabla f)(\sigma(u))\ge 0$. This shows that $|\nabla
f|^2(\sigma(u))=|\nabla f|^2(\sigma(0))$ for $u\le 0$ since
$|\nabla f|^2$ has its minimum at $o$. Hence we have that
$S(\sigma(u))=\max_M S$ for all $u\le 0$. However since
$-f(\sigma(u))=-f(\sigma(0))+\int_u^0 |\nabla f|^2\,
du=-f(\sigma(0))-u|\nabla f|^2(\sigma(0))\to +\infty$ as $u\to
-\infty$ we can conclude that $\sigma(u)\to \infty$. This is a
contradiction with the assumption that $\limsup_{x\to \infty}S(x)
<\max_M S$. Hence we have that $\lambda=\max_M S$ which implies
$\nabla f=0$ at $o$ and
$$
\liminf_{x\to \infty} |\nabla f|^2 \ge 2\eta^2\doteqdot
\lambda-\limsup_{x\to \infty} S(x)>0.
$$
By considering any minimizing geodesic $\gamma(s)$ emitting from $o$
and the fact $\frac{d}{ds}(f(\gamma(0)))=0$ and
$\frac{d^2}{ds^2}(f(\gamma(s)))\ge 0$, it is clear that $o$ is the
minimal point of $f$ and $\langle \nabla f, \nabla r\rangle (x)\ge
0$ for any $x\in M\setminus \{o\}$. Let $R_0$ be such that $|\nabla
f|^2(x)\ge \eta^2$ for all $x\in M\setminus B(o, R_0)$. Consider
again an integral curve $\sigma(u)$ passing $x$. Since $|\nabla f|$
is bounded and $M$ is complete, the curve is defined for all
$-\infty < u< +\infty$. Notice that $\sigma(u)\in B(o, r(x))$ for
all $u\le 0$ and
$$
f(\sigma(0))-f(\sigma(u))=\int_u^0 |\nabla f|^2\, du \ge (-u)\eta^2
$$
as along as $\sigma(u)\in M\setminus B(o, R_0)$. From this we infer
that there exist some $u_0$ such that $\sigma(u_0)\in B(o, R_0)$. On
the other hand
$$
f(x) = f(\sigma(u_0))+\int_{u_0}^0 |\nabla f|^2\, du \ge
f(\sigma(u_0))+\eta \int_{u_0}^0 |\sigma'(u)|\, du \ge
f(\sigma(u_0))+\eta d(x, \sigma(u_0)).
$$
This implies the desired lower estimate. The final conclusion
follows easily from the above estimate on $|\nabla f|$ and the
convexity of $f$. \qed

\begin{remark}
If the sectional curvature of $(M, g)$ is nonnegative, one can
show that the claim of the lemma  holds under the assumption that
$S(o)=\max_M S$, as far as $M$ does not admit any flat factor
$\R^k$. The reason is the following. First if the claimed result
fails, one can conclude that  $f_{ij}$ has an eigenvector
corresponding to the zero eigenvalue somewhere. Note that for the
associated Ricci flow, the function  $f(x, \tau)$, defined as the
pull back via the diffeomorphism generated by $\nabla f$,
satisfies the heat equation (cf. {\rm \cite{CLN}} for details).
Then the result follows from the strong tensor maximum principle
and a splitting theorem of noncompact manifolds proved in {\rm
\cite{N0}}.
\end{remark}

In the both cases $o$ is a minimum point of $f$ and $\lambda$ is a
geometric invariant, namely $\max_{x\in M} S(x)$. Also we have
seen that both $|\nabla f|$ and $|\Delta f|$ are bounded. We
normalize $f$ so that $\int_M e^{-f}\, d\Gamma=1$. Integration by
parts gives the following weighted Poincar\'e inequality.

\begin{proposition} Let $(M, g, f)$ be a gradient steady soliton. Then for
any compact supported smooth function $u=e^{-\psi}$ with $\int_M
u\, d\Gamma=1$, we have that
$$
\int_M \left(|\nabla \psi|^2-3S\right)u\, d\Gamma \ge -\lambda.
$$
\end{proposition}
\noindent {\sl Proof.}
 The proof follows from the following simple calculation:
\begin{eqnarray*}
\int_M \left(|\nabla \psi|^2-3S\right)u\,  &=& \int_M \left(|\nabla \psi|^2-2\langle \nabla \psi, \nabla f\rangle +
|\nabla f|^2+2\Delta f -|\nabla f|^2-3S\right)u\\
&\ge& \int_M \left(2\Delta f -|\nabla f|^2-3S\right)u=-\lambda,
\end{eqnarray*}
for all normalized $u$.\qed

This is a sharp inequality, at least it is so under the assumption of Lemma
\ref{gs-help1}, since for this case the equality holds when $u=e^{-f}$. An
equivalent form is that
$$
\int_M \left(4|\nabla \varphi|^2 -3S \varphi^2\right) \, d\Gamma
\ge -\lambda \int_M \varphi^2\, d\Gamma
$$
for any $\varphi\in L^2(M)$. The weighted Poincar\'e inequality
and its geometric meanings have recently been studied in
\cite{LW}.

\section{Extensions and an analogue of the c-theorem}

For the re-normalization group flow, there exists the so-called
central charge $c(t)$ invariant \cite{Zam} for the flow  such that
{\it $c(t)$ is monotone non-increasing in $t$}. Moreover {\it $c(t)$ is
always nonnegative}. For Ricci flow, there are Perelman's monotonic
quantities such as the $\mathcal{W}(g, \sigma, f)$-entropy,
defined as
$$
\mathcal{W}(g, \sigma, \varphi)\doteqdot \int_M \left(\sigma
(|\nabla \varphi|^2+S)+\varphi-n\right) u\, d\Gamma
$$
for any $u=e^{-\varphi}/{(4\pi \sigma)^{\frac{n}{2}}}$ with
$\int_M u=1$, and  associated $\mu(g, \sigma)\doteqdot
\inf_{\int_M u=1}\mathcal{W}(g, \sigma, \varphi)$,
$\nu(g)\doteqdot \inf_{\sigma>0} \mu(g, \sigma)$ invariants, as
well as  the so-called reduced volume. The quantity
$\mathcal{W}(g, \sigma, \varphi) $, $\mu(g, \sigma)$ and $\nu(g)$
may not be finite when $M$ is not compact. Proposition
\ref{as-lsi1} and Proposition \ref{as-lsi2} ensures that is the
case for the shrinkers with bounded nonnegative Ricci curvature.
The reduced volume \cite{P} is always nonnegative by the
definition. But it is monotone non-decreasing instead of
non-increasing in $t$ (along the flow). Utilizing the sharp LSI's
proved for the shrinkers and expanders we shall show in this
section that the logarithmic Sobolev constants $\mu_s$ and $\mu_e$
are nonnegative, at least for the gradient shrinking/expanding
solitons (with some mild assumptions on the Ricci curvature). In
view of the monotonicity of the entropy, and the fact that  the
gradient shrinking solitons often arises at the singularity, one
can  view the monotonicity of the entropy together with  the
result proved here as  an analogue of the c-theorem. Namely, for
the solution to the Ricci flow, one can view $-\mu(g, \tau)$ as
the analogue of the $c(t)$-invariant. Perelman's entropy formula
concludes that it is monotone non-increasing. Our result concludes
that $-\mu(g, 1)=\mu_s$ and it is nonnegative. One should note
that if the $\nu(g(t))$ invariant of Perelman \cite[Section 3]{P}
is well-defined/finite, unfortunately this is not always the case,
then at least for the compact manifolds, the $-\nu(g(t))$ would be
nonnegative and non-increasing along the Ricci flow.

We shall show two results on the sign of the invariants $\mu_s$
and $\mu_e$. The case of $\mu_e$ is an easy application of a
rigidity result in \cite{N1}.

\begin{proposition} Let $(M, g, f)$ be a gradient expanding soliton with $\Ric \ge 0$.
Then  $\mu_e\ge 0$. If   $\mu_e=0$  then $(M, g)$ must be isometric to $\R^n$.
\end{proposition}
\noindent {\sl Proof.} Assuming that $\mu_e\le 0$,  Theorem
\ref{ges-main1} then implies that
$$
\int_M \left(|\nabla \psi|^2+\psi-n\right)\rho\, d\Gamma \ge 0.
$$
Then by \cite[Theorem 1.4]{N1}, one can see a detailed account in
\cite[pages 314--333]{Chowetc}, we can conclude that $(M, g)$ is
isometric to $\R^n$, on which $\mu_e=0$. \qed

Similar result holds for gradient  shrinking solitons. For that we
have to assume that the  curvature tensor of $(M, g)$ is
bounded.

\begin{theorem}\label{ctheorem}
Let $(M, g)$ be a gradient shrinking soliton with bounded
curvature. Let $f$ be the normalized potential function as before,
then $\mu_s\ge 0$.
\end{theorem}

\begin{remark}
After the appearance of our paper, Yokota  \cite{Ta} generalized
the above result by assuming only the lower bound of Ricci
curvature. The proof makes uses of Perelman's reduced volume.
\end{remark}

\noindent {\sl Proof.} We first prove the result under the extra
assumption that $\Ric\ge 0$. Recall from the introduction that
there is an associated  solution $g(t)$ (with $-\infty < t<0$,
$t=\eta-1$) to Ricci flow generated by pulling back the metric via
the diffeomorphisms generated by the vector field $\nabla f$. The
original metric $g$ corresponds to the one $g(-1)$ (meaning
$t=-1$). Proposition \ref{as-lsi1} and Proposition \ref{as-lsi2}
imply that $\mu(g(-1), \sigma)$ and $\nu(g(-1))$ are finite. Since
$g(t)$ is just the re-scale of $g(-1)$, we have that for any
$-\infty<t<0$, $\mu(g(t), \sigma)$, $\nu(g(t))$ are also finite.
Now let $H(y, t; x, t_0)$ (with $t<t_0<0$) be the (minimal)
positive fundamental solution to the conjugate heat equation:
$$
\left(-\frac{\partial}{\partial t}-\Delta_y +S(y, t)\right)H(y, t; x, t_0)=0
$$
being the $\delta_x(y)$ at $t=t_0$. By a result of Perelman
\cite[Corollary 9.3]{P}, see also \cite{CTY,N-LYH}, we know that
$$
v_H(y, t) \doteqdot (t_0-t)\left(2\Delta \varphi -|\nabla \varphi|^2+S\right)+\varphi-n\le 0
$$
with $H(y, t; x, t_0)=e^{-\varphi(y, t)}/{(4\pi
(t_0-t))^{\frac{n}{2}}}$.  This  implies in particular
$$
\mu(g(-1), t_0+1)\le \int_M v_H(y, -1) H(y, -1)\,
d\Gamma_{g(-1)}\le   0.
$$
Here to ensure  the inequality $v_H\le 0$  the extra assumption
that the curvature tensor of $M$ is uniformly bounded  is needed
\cite{CTY}.

On the other hand Theorem \ref{main1} asserts that $\mu(g(-1),
1)\ge -\mu_s$. The result would follow if we show that $\mu(g(-1),
t_0+1)\to \mu(g(-1), 1)$ as $t_0\to 0$. For $t_{0, i}\to 0$,
consider  minimizers $\varphi_i$ of $\mathcal{W}(g(-1),
1+t_{0, i}, \varphi)$ (for simplicity we write $g(-1)$ back to $g$
from now on). Let $\sigma_i=1+t_{0, i}\to 1$. We assume that
$\frac{1}{2}\le \sigma_i \le 1$. By Proposition \ref{as-lsi2}  and
the above we have that
$$
0\ge \mu(g,\sigma_i)\ge -\mu_s-nA.
$$
Write $w_i=e^{-\varphi_i/2}$. Then, essentially from definition, the $w_i\in W^{1,2}(M)$. The Euler-Lagrangian equation is
\begin{equation}\label{el-i}
-4\sigma_i \Delta w_i +\sigma_i S w_i -nw_i -2w_i \log w_i =\mu(g, \sigma_i)w_i
\end{equation}
for $\int_M w_i^2=(4\pi \sigma_i)^{\frac{n}{2}}\le (4\pi)^{\frac{n}{2}}$.
Integrating over $M$ we have that
$$
4\sigma_i\int_M |\nabla w_i|^2 =\mu(g, \sigma_i)(4\pi \sigma_i)^{\frac{n}{2}}+\int_M \left(w_i^2\log w_i^2+nw^2_i-\sigma_i S w_i^2\right)
$$
which implies
\begin{equation}\label{grad1}
4\sigma_i\int_M |\nabla w_i|^2 \le \int_M w_i^2\log w_i^2 +n(4\pi)^{\frac{n}{2}}
\end{equation}

On the other hand, writing ${w_i^2}/{(4\pi
\sigma_i)^{\frac{n}{2}}}={e^{-\tpsi}}/{\pi^{\frac{n}{2}}}$ and
using that $ \mathcal{W}(g, \frac{1}{4}, \tpsi)\ge \mu(g,
\frac{1}{4})$,
$$
\int_M |\nabla w_i|^2 \ge (4\pi\sigma_i)^{\frac{n}{2}}\mu(g, \frac{1}{4})-\frac{nA}{4}(4\pi \sigma_i)^{\frac{n}{2}}
+\int_M w_i^2\log w_i^2-\frac{n}{2}\log (4\pi \sigma_i). $$
Combining with (\ref{grad1}), one can find $C=C(A,n)$ such that
$$
\int_M |\nabla w_i|^2 \le C(A, n)
$$
which implies that $\|w_i\|_{W^{1,2}(M)}$ is uniformly bounded. It then
 implies that $w_i\to w_\infty$ in the the dual norm of
$W^{1,2}(M)$ and strongly in $L^2(M)$, for some $w_\infty\in
W^{1,2}(M)$. Due to the bound $\mu(g, \sigma_i)$ we may also
assume that $\mu(g, \sigma_i)\to \mu_{\infty}(g)$. Clearly $\mu_\infty(g)\le 0$. It is evident
that
 $\int_M w_\infty^2=(4\pi)^{\frac{n}{2}}$. We shall show that on every compact subset $K$, after passing to subsequences,  $w_i$ converges to  $w_\infty$, say in $C^0$-fashion. This will imply that $w_\infty$ satisfies the equation
$$
-4 \Delta w_\infty + S w_\infty -nw_\infty -2w_\infty \log w_\infty =\mu_\infty(g) w_\infty.
$$
Integration by parts yields that $$
\int_M \left(4|\nabla w_\infty|^2 +S w_\infty^2 -2w_\infty^2\log w_\infty-nw_\infty^2\right) =\mu_\infty(g)(4\pi)^{\frac{n}{2}}.
$$
This implies that $\mu_\infty(g) \ge \mu(g, 1)$, which is
enough to conclude that $\mu_s\ge 0$ since $0\ge
\mu_{\infty}(g)\ge \mu(g, 1) = -\mu_s$. The claim that $w_i\to
w_\infty$ in $C^0$ norm can be proved using Sobolev embedding
theorem (over compact region $K$), interior $L^p$-estimates, and the
compactness of the Sobolev embedding. Since it is rather standard we leave the details to the
interested reader. One can also find this in the forthcoming book
\cite{Chow2}.

Now we point out how one can modify the above argument to the
general case. In fact in the proof above the assumption that
$\Ric\ge 0$ is only used, via Proposition \ref{as-lsi2}, to ensure
that $\mu(g, \sigma)$ is uniformly bounded for $1-\delta \le
\sigma <1$, for some $\delta>0$. This can be done for the case
that $|\Ric|\le A$ for some $A>0$. We state this as a separate
result below. \qed

\begin{proposition}
Assume that on a complete Riemannian manifold $(M, g)$,  $\mu(g,
1)>-\infty$ and $\Ric\ge -A$ and $S\le B$ for some positive
numbers $A$ and $B$. Then for any $0<\sigma<1$,
\begin{equation}
\mu(g, \sigma) \ge \mu(g, 1) -nA\sigma-B -\left( \frac{ A^2n}{2}+An\right)(1-\sigma).
\end{equation}
\end{proposition}
\noindent {\sl Proof.} As in Proposition \ref{as-lsi2}, $\mu_0(g,
1)\ge \mu(g, 1) -B$. Let $u_0(x)={e^{-\tilde{\psi}}}/{(4\pi
\sigma)^{n/2}}$ be a smooth function with compact support such
that $\int_M u_0 =1$. Similarly let $u(x, t)={e^{-\varphi}}/{(4\pi
\tau)^{n/2}}$ be the solution to the heat equation with $u(x,
0)=u_0(x)$. Here $\tau(t)=\sigma+t$. We shall use the entropy
formula from \cite{N1} to estimate
$$
\mathcal{W}_0(0)\doteqdot \int_M \left(\sigma |\nabla \tilde \psi|^2 +\tilde \psi -n\right)u_0.
$$
Let $F(t)=\int_M |\nabla \varphi|^2 u$. The entropy formula of
\cite{N1} implies that the entropy
$$
\mathcal{W}_0(t)\doteqdot \int_M \left(\tau |\nabla \varphi|^2+\varphi -n\right) u
$$
satisfies the estimate
\begin{align*}
\frac{d} {d t}  \mathcal{W}_0(t) &\le   -2\tau \int_M \left|\nabla_i\nabla_j \varphi -\frac{1}{2\tau}g_{ij}\right|^2 u +2\tau A F(t) \le -\frac{2\tau}{n} \int_M \left(\Delta \varphi -\frac{n}{2\tau}\right)^2 u +2\tau A F(t)\\
&\le -\frac{2\tau}{n} \left(\int_M (\Delta \varphi
-\frac{n}{2\tau})u\right)^2+2\tau A F(t) =
-\frac{2\tau}{n}\left(F(t)-\frac{n}{2\tau}\right)^2 +2\tau A F(t).
\end{align*}
Viewing  the right hand side above as a quadratic polynomial in $X=F(t)-\frac{n}{2\tau}$, by an elementary consideration we have that
$$
\frac{d}{dt}\mathcal{W}_0(t)\le \frac{ A^2n}{2}+nA
$$
for $\tau \le 1$.
Hence
$$
\mathcal{W}_0(0)\ge \mathcal{W}_0(1-\sigma) -\left( \frac{ A^2n}{2}+An\right)(1-\sigma).
$$
This shows that
$$
\mu_0(g, \sigma) \ge \mu_0(g, 1)-\left( \frac{ A^2n}{2}+An\right)(1-\sigma).
$$
Finally we have that $\mu(g, \sigma)\ge \mu(g, 1)-nA\sigma-B -\left( \frac{ A^2n}{2}+An\right)(1-\sigma).$
\qed

When $f=\operatorname{constant}$,  $(M, g)$ is a compact Einstein manifold
with $\Ric_M =\frac{1}{2}g_M$. The theorem concludes that
$$
\mu_s=\frac{n}{2}-\log (V(M))+\frac{n}{2}\log (4\pi) \ge 0
$$
where $V(M)$ is the volume of $(M, g_M)$. Among all such manifolds
the sphere $\Sph^n$ has the smallest $\mu_s$. In this case $\mu_s$
is monotone non-increasing in $n$ and has the  limit
$\frac{1}{2}\log \frac{e}{2}$ as $n\to \infty$, at least for the
case that $n$ is even. In fact,
$$
\mu_s(\Sph^{2k})=\log \frac{e^k (2k-1)!}{\left(2(2k-1)\right)^k
(k-1)!}.
$$
It is also easy to see that $\mu_s(\R^n)=0$ and $\mu_s(M_1\times
M_2)=\mu_s(M_1)+\mu_s(M_2)$.

\begin{remark}
If $(M, g(t))$ is a solution to Ricci flow on compact manifold $M$
over $[0, T)$. Then for any $0\le t_1<T$, $\mu(g(t_1), T-t_1)\le
0$ by an argument similar as (but easier than) the above.  For the
steady gradient soliton, it is clear that $\lambda\ge 0$ for any
steady solitons with $S\ge 0$. We conjecture that if $\mu_s=0$,
then the shrinker has to be isometric to $\R^n$. In \cite{Ta},
this conjecture has been proved recently.
\end{remark}

\begin{corollary} Let $(M, g, f)$ be a gradient shrinking soliton as in Theorem \ref{ctheorem}. Then
$$
\int_M f\frac{e^{-f}}{(4\pi)^{n/2}}\le \frac{n}{2}.
$$
\end{corollary}

\section*{Acknowledgments} {
JAC acknowledges partial support from the project
MTM2008-06349-C03-03 from DGI-MCINN (Spain) and IPAM (UCLA) where
this work was essentially done. LN was supported in part by NSF
grant DMS-0805834, Institut Henri Poincar\'{e} and an Alfred P.
Sloan Fellowship, USA. LN would like to thank Dan Friedan for the
lecture \cite{Fri} on the re-normalization group flow and
information regarding the c-theorem. This motivated  Theorem
\ref{ctheorem}. }

\bibliographystyle{amsalpha}

{\sc Addresses}

{\sc Jos\'e A. Carrillo}, ICREA and Departament de Matem\`atiques,
Universitat Aut\`onoma de Barcelona, E-08193-Bellaterra, Spain

email: carrillo@mat.uab.es

{\sc Lei Ni},
 Department of Mathematics, University of California at San Diego, La Jolla, CA 92093, USA


email: lni@math.ucsd.edu

\end{document}